\documentclass[a4paper,12pt,reqno]{amsart}
\usepackage{amsmath}
\usepackage{amssymb}
\usepackage{amsthm}

\usepackage{amscd}
\usepackage{amsfonts}
\usepackage{setspace}
\usepackage{version}

\title[Riemann zeta function and log-correlated random variables]{A note on the maximum of the Riemann zeta function, and log-correlated random variables}
\author{Adam J Harper}
\address{Centre de recherches math\'{e}matiques, Universit\'{e} de Montr\'{e}al, Pavillon Andr\'{e}-Aisenstadt, 2920 Chemin de la tour, bur. 5357, Montr\'{e}al QC H3T 1J4, Canada}
\email{harperad@crm.umontreal.ca}
\date{2nd April 2013}
\thanks{The author is supported by a postdoctoral fellowship from the Centre de recherches math\'{e}matiques in Montr\'{e}al.}

\numberwithin{equation}{section}
\theoremstyle{plain}

\onehalfspace
\newcommand{\R}{\mathbb{R}}
\newcommand{\E}{\mathbb{E}}
\newcommand{\p}{\mathbb{P}}

\newcommand{\C}{\mathbb{C}}

\newtheorem{conj1}{Conjecture}
\newtheorem{prop1}{Proposition}
\newtheorem{prop2}[prop1]{Proposition}
\newtheorem{conj2}[conj1]{Conjecture}
\newtheorem{tail}{Tail Inequality}
\newtheorem{lem1}{Lemma}
\newtheorem{lem2}[lem1]{Lemma}
\newtheorem{lowerb}{Lower Bound}
\newtheorem{lem3}[lem1]{Lemma}
\newtheorem{comparison}{Comparison Inequality}
\newtheorem{clt}{Central Limit Theorem}

\begin{document}

\maketitle

\begin{abstract}
In recent work, Fyodorov and Keating conjectured the maximum size of $|\zeta(1/2+it)|$ in a typical interval of length $O(1)$ on the critical line. They did this by modelling the zeta function by the characteristic polynomial of a random matrix; relating the random matrix problem to another problem from statistical mechanics; and applying a heuristic analysis of that problem.

In this note we recover a conjecture like that of Fyodorov and Keating, but using a different model for $|\zeta(1/2+it)|$ in terms of a random Euler product. In this case the probabilistic model reduces to studying the supremum of Gaussian random variables with logarithmic correlations, and can be analysed rigorously.
\end{abstract}

\section{Introduction}
When $\Re(s) > 1$, the Riemann zeta function can be expressed as an absolutely convergent Euler product:
$$ \zeta(s) = \prod_{p \; \textrm{prime}} \left(1- \frac{1}{p^{s}}\right)^{-1} . $$
Thus $\zeta(s) \neq 0$ when $\Re(s) > 1$. For other $s$ we do not have such a nice expression for $\zeta(s)$, and its behaviour remains quite mysterious. However, one can still approximate $\zeta(s)$ in useful ways. Gonek, Hughes and Keating~\cite{gonekhugheskeating} showed, very roughly speaking, that if the Riemann Hypothesis is true then
$$ \zeta(1/2+it) \approx \prod_{p \leq X} \left(1- \frac{1}{p^{1/2+it}}\right)^{-1} \cdot \prod_{|\gamma - t| < 1/\log X} (ci(t-\gamma)\log X), $$
on a wide range of the parameter $X$ (relative to $t$). Here $c > 0$ is a constant, and the second product is over ordinates $\gamma$ of zeros of the zeta function.

The ordinates $\gamma$ and the numbers $p^{it}$ are difficult to analyse rigorously, but it is widely believed that they behave like certain random objects. One might imagine that the $p^{it}$ behave, for most $t$, like independent random variables $U_{p}$ distributed uniformly on the unit circle, and Selberg's central limit theorem for $\log|\zeta(1/2+it)|$ provides rigorous support for that belief. The correct model for the $\gamma$ is believed to be the eigenvalues of a random unitary matrix. Inserting these random objects in the ``hybrid Euler--Hadamard product'' above, one obtains (for any choice of $X$) a random model for $\zeta(1/2+it)$.  Farmer, Gonek and Hughes~\cite{farmergonekhughes} analysed the likely (i.e. with probability $1-o(1)$) behaviour of the maximum of $T$ independent copies of that random model, and so were led to a conjecture about $ \max_{0 \leq t \leq T} |\zeta(1/2+it)|$. It turns out that one derives the same conjecture for a wide range of choices of $X$.

More recently, Fyodorov and Keating~\cite{fyodkeat} considered a short interval version of the maximum, namely $\max_{T \leq t \leq T + 2\pi} |\zeta(1/2+it)|$. See also the announcement~\cite{fyodhiarykeat} by Fyodorov, Hiary and Keating. Amongst other things, they make a conjecture that may be interpreted\footnote{Fyodorov and Keating~\cite{fyodkeat} actually make an even more precise conjecture, about the distribution of $\log\max_{T \leq t \leq T + 2\pi} |\zeta(1/2+it)| - \log\log T + (3/4)\log\log\log T$ as $T$ varies.} as follows:
\begin{conj1}[Fyodorov and Keating, 2012]
Let $\epsilon > 0$, and suppose that $T_{1}$ is large enough depending on $\epsilon$. Then for all $T_{1} \leq T \leq 2T_{1}$, except for a set of ``bad'' $T$ of measure at most $\epsilon T_{1}$, we have
$$ \max_{T \leq t \leq T + 2\pi} |\zeta(1/2+it)| = \exp(\log\log T - (3/4)\log\log\log T + O_{\epsilon}(1)) , $$
where $O_{\epsilon}(1)$ denotes a quantity bounded in terms of $\epsilon$.
\end{conj1}
As well as being interesting in itself, Fyodorov and Keating~\cite{fyodkeat} observe that their conjecture is easier to investigate numerically than the conjecture about $ \max_{0 \leq t \leq T} |\zeta(1/2+it)|$, because one only needs to calculate the behaviour of the zeta function in several randomly chosen intervals of length $2\pi$, rather than in a complete interval of length $T$. See Fyodorov and Keating's paper for a discussion of such numerical experiments.

As in the work of Farmer, Gonek and Hughes~\cite{farmergonekhughes} on $\max_{0 \leq t \leq T} |\zeta(1/2+it)|$, and in huge amounts of previous work on e.g. the moments $\int_{0}^{T} |\zeta(1/2+it)|^{2k} dt$ of the zeta function, Fyodorov and Keating are led to their conjectures by modelling $|\zeta(1/2+it)|$ using the characteristic polynomials of random matrices. More precisely, they speculate that $\max_{T \leq t \leq T + 2\pi} |\zeta(1/2+it)|$ will typically behave in the same way as
$$ \max_{0 \leq \theta \leq 2\pi} |p_{N}(\theta)| , $$
where $p_{N}(\theta)$ is the characteristic polynomial of a random $N \times N$ unitary matrix, and $N$ is the integer closest to $\log T$. However, it seems difficult to precisely analyse the random variable $\max_{0 \leq \theta \leq 2\pi} |p_{N}(\theta)|$ in a rigorous way. As Fyodorov and Keating explain, when Farmer, Gonek and Hughes~\cite{farmergonekhughes} encounter $\max_{0 \leq \theta \leq 2\pi} |p_{N}(\theta)|$ they only require some information about its extreme tail behaviour, because they then take the maximum of $T$ independent copies. Fyodorov and Keating~\cite{fyodkeat} require precise distributional information, and end up relying on a (very interesting) heuristic analysis based on comparing $\max_{0 \leq \theta \leq 2\pi} |p_{N}(\theta)|$ with some statistical mechanics problems.

\vspace{12pt}
In this note we investigate $\max_{T \leq t \leq T + 2\pi} |\zeta(1/2+it)|$ using a random Euler product model, that roughly corresponds to choosing $X$ as a power of $T$ in the hybrid Euler--Hadamard product. The conjecture we derive turns out to be essentially the same as Conjecture 1, and the probabilistic analysis of this model can be performed rigorously, which hopefully serves as additional evidence in favour of Conjecture 1.

First we show that, except on a set of small measure (essentially consisting of points very close to zeta zeros), $\log|\zeta(1/2+it)|$ is very close to $\Re(\sum_{p \leq T} \frac{1}{p^{1/2+it}} \frac{\log(T/p)}{\log T})$. We do this by adapting an argument of Soundararajan~\cite{soundmoments} that gave an upper bound for $\log|\zeta(1/2+it)|$, which is valid for all $t$ and which we will also need. These results assume the Riemann Hypothesis, but this seems acceptable in our heuristic context.
\begin{prop1}[Adapted from the Main Proposition of Soundararajan~\cite{soundmoments}]
Assume the Riemann Hypothesis is true, and let $T$ be large. Then for any $T \leq t \leq 2T$ we have
$$ \log|\zeta(1/2+it)| \leq \Re(\sum_{p \leq T} \frac{1}{p^{1/2 + 1/\log T + it}} \frac{\log(T/p)}{\log T} + \sum_{p^{2} \leq T} \frac{(1/2)}{p^{1 + 2/\log T + 2it}} \frac{\log(T/p^{2})}{\log T}) + O(1) . $$

Moreover, there exists a set $\mathcal{H} \subseteq [T,T+2\pi]$, of measure at least $1.99 \pi$, such that
$$ \log|\zeta(1/2+it)| = \Re(\sum_{p \leq T} \frac{1}{p^{1/2+it}} \frac{\log(T/p)}{\log T}) + O(1) \;\;\;\;\; \forall t \in \mathcal{H} . $$
\end{prop1}
We will prove Proposition 1 in $\S 2$.

The other part of our analysis is the following probabilistic result, which we shall discuss further in just a moment.
\begin{prop2}
Let $T$ be large. Let $(U_{p})_{p \leq T}$ be independent random variables, each distributed uniformly on the unit circle in $\C$. Then
$$ \p\biggl(\max_{0 \leq h \leq 2\pi} \Re(\sum_{p \leq T} \frac{U_{p}}{p^{1/2+1/\log T +ih}} \frac{\log(T/p)}{\log T} + \frac{1}{2} \sum_{p^{2} \leq T} \frac{U_{p}^{2}}{p^{1+2/\log T +2ih}} \frac{\log(T/p^{2})}{\log T}  ) $$
$$ \leq \log\log T - (1/4)\log\log\log T + O(\sqrt{\log\log\log T}) \biggr) $$
is $1 - o(1)$ as $T \rightarrow \infty$.

Moreover, if $\mathcal{H} \subseteq [0,2\pi]$ is any fixed set of measure at least $1.99 \pi$ (say) then
$$ \p(\max_{h \in \mathcal{H}} \Re(\sum_{p \leq T} \frac{U_{p}}{p^{1/2+ih}} \frac{\log(T/p)}{\log T} ) \geq \log\log T - 2\log\log\log T - O((\log\log\log T)^{3/4})) = 1 - o(1) . $$
\end{prop2}

As discussed previously, it seems reasonable to assume that for a ``typical'' value of $T$ the set of values $(p^{iT})_{p \leq T}$ will behave, in an average sense, like $(U_{p})_{p \leq T}$. Thus it seems reasonable to assume that, for a typical value of $T$, $\left(\Re(\sum_{p \leq T} \frac{1}{p^{1/2+it}} \frac{\log(T/p)}{\log T}) \right)_{T \leq t \leq T + 2\pi}$ will behave in the same way as $\left(\Re(\sum_{p \leq T} \frac{U_{p}}{p^{1/2+ih}} \frac{\log(T/p)}{\log T}) \right)_{0 \leq h \leq 2\pi}$. Combining this assumption with Propositions 1 and 2, we are led to the following conjecture:
\begin{conj2}
Let $\epsilon > 0$, and suppose that $T_{1}$ is large enough depending on $\epsilon$. Then for all $T_{1} \leq T \leq 2T_{1}$, except for a set of ``bad'' $T$ of measure at most $\epsilon T_{1}$, we have
\begin{eqnarray}
&& \exp(\log\log T - (2 + o(1))\log\log\log T) \nonumber \\
& \leq & \max_{T \leq t \leq T + 2\pi} |\zeta(1/2+it)| \leq \exp(\log\log T - (1/4+ o(1))\log\log\log T) \nonumber
\end{eqnarray}
\end{conj2}
Thus, as claimed, we were led to essentially the same conclusion\footnote{It is known that $\log|p_{N}(\theta)|$ can be expressed as a trigonometric series in $\theta$, whose coefficients are random variables any finite number of which have independent Gaussian limiting distributions (as $N \rightarrow \infty$). See Diaconis and Shahshahani's paper~\cite{diaconisshah}, and also the discussion in Fyodorov and Keating's paper~\cite{fyodkeat}. Thus it may not be too surprising that Conjecture 2, which is derived from studying random trigonometric sums, matches Conjecture 1, which reflects the presumed behaviour of $\log|p_{N}(\theta)|$.} as Conjecture 1. If we could prove a more precise version of Proposition 2 then a more precise version of Conjecture 2 would follow.

\vspace{12pt}
We conclude this introduction with some remarks about Proposition 2. If $|h_{1}-h_{2}| \leq 1/\log T$ then the random variables
$$ X(h_{1}) := \Re(\sum_{p \leq T} \frac{U_{p}}{p^{1/2+ih_{1}}} \frac{\log(T/p)}{\log T} ), \;\;\;\;\; X(h_{2}) := \Re(\sum_{p \leq T} \frac{U_{p}}{p^{1/2+ih_{2}}} \frac{\log(T/p)}{\log T} ), $$
are almost perfectly correlated, since $p^{ih_1} \approx p^{ih_2}$ for all $p \leq T$. When $0 \leq h_{1}, h_{2} \leq 2\pi$ are further apart we have
\begin{eqnarray}
\E X(h_{1}) X(h_{2}) & = & \sum_{p_{1}, p_{2} \leq T} \frac{\E \Re(U_{p_{1}} p_{1}^{-ih_{1}}) \Re(U_{p_{2}} p_{2}^{-ih_{2}}) }{p_{1}^{1/2} p_{2}^{1/2}} \frac{\log(T/p_{1}) \log(T/p_{2})}{\log^{2}T} \nonumber \\
& = & \frac{1}{2} \sum_{p \leq T} \frac{\cos((h_{1}-h_{2})\log p)}{p} \frac{\log^{2}(T/p)}{\log^{2}T} \nonumber \\
& \approx & (1/2)\log(1/|h_{1}-h_{2}|) , \nonumber
\end{eqnarray}
where the second equality follows by writing $\Re(U_{p} p^{-ih}) = (1/2)(U_{p}p^{-ih} + \overline{U_{p}}p^{ih})$ (and noting that $\E U_{p_{1}} U_{p_{2}} = 0$ for all $p_{1},p_{2}$), and the third line from standard estimates for sums over primes. See the first appendix, below, for a precise calculation. Unsurprisingly given the origin of the $X(h)$, this logarithmic type of covariance structure matches the two-point correlation function of $\log|\zeta(1/2+i(T+h))|$, as discussed in Fyodorov and Keating's paper~\cite{fyodkeat}. Such covariance structures also appear in connection with branching random walks, with the Gaussian Free Field, and in general with random variables indexed by trees. See Zeitouni's lecture notes~\cite{zeitouni} for some discussion of the former problems, and also see the references cited by Fyodorov and Keating.

The upper bound in Proposition 2 is, in principle, quite straightforward to obtain, but complications arise because the maximum is taken over an interval rather than a discrete set, and in trying to obtain the subtracted term $-(1/4)\log\log\log T$. We will prove it in $\S 3.1$. The lower bound is harder, but one can use a quantitative form of the multivariate central limit theorem to replace the random variables $\left(\Re(\sum_{p \leq T} \frac{U_{p}}{p^{1/2+ih}} \frac{\log(T/p)}{\log T}) \right)_{0 \leq h \leq 2\pi}$ by Gaussian random variables with the same means and covariances, and then use a general lower bound inequality from the author's paper~\cite{harpergp}. See $\S 3.3$, below. In fact that paper studies random variables very much like the $X(h)$, in connection with a problem of Hal\'{a}sz on random multiplicative functions. (In \cite{harpergp} the $U_{p}$ are replaced by real random variables, and the range of $h$ is a little different, but otherwise the situations are very similar.) The author hopes that the exposition here will be a useful supplement to that in \cite{harpergp}, where most of the focus was on establishing the basic inequality rather than the application.

As Fyodorov and Keating~\cite{fyodkeat} briefly discuss, (and see $\S 3.1$, below), if the values $\log|\zeta(1/2+i(T+h))|$ behaved ``independently'' when $|h_1 - h_2| \gg 1/\log T$ then the subtracted term $-(3/4)\log\log\log T$ in Conjecture 1 would be incorrect, and instead one would have a subtracted term $-(1/4)\log\log\log T$. The term $-(3/4)\log\log\log T$ is believed to be a manifestation of the logarithmic correlations of the $\log|\zeta(1/2+i(T+h))|$ (and the corresponding random models). Thus it seems an interesting problem to sharpen Proposition 2, and to rigorously analyse the random matrix model that motivated Fyodorov and Keating. It is also conceivable that one could obtain rigorous results, in the direction of Conjecture 1, about the zeta function itself. For example, it is well known that if $T \leq t \leq 2T$ then $\zeta(1/2+it) = \sum_{n \leq T} 1/n^{1/2 +it} + O(1/\sqrt{T})$, and so
$$ \sum_{T_1 \leq n \leq 2T_1} \max_{n \leq t \leq n+2\pi} |\zeta(1/2+it)|^{2} \ll \sum_{T_1 \leq n \leq 2T_1} \max_{n \leq t \leq n+2\pi} \left|\sum_{n \leq T_{1}} \frac{1}{n^{1/2+it}} \right|^{2} + O(1) \ll T_{1} \log^{2}T_{1} , $$
using a discrete mean value result for Dirichlet polynomials (as in e.g. Theorem 5.3 of Ivi\'{c}~\cite{ivic}). This implies that $\max_{T \leq t \leq T+2\pi} |\zeta(1/2+it)| = O_{\epsilon}(\log T)$, except for a set of bad $T_1 \leq T \leq 2T_1$ of measure at most $\epsilon T_1$. However, it is not clear how to obtain any information at the level of the $\log\log\log T$ corrections in Conjectures 1 and 2.

\section{The number theoretic part : Proof of Proposition 1}
The first part of Proposition 1, giving an upper bound for $\log|\zeta(1/2+it)|$, is a special case of the Main Proposition of Soundararajan~\cite{soundmoments} (choosing $x=T$ and $\lambda = 1$ there). Thus we will just prove the second part of the proposition.

Let $C, K > 0$ be large absolute constants, suitable values of which could be extracted from the following arguments. We will begin by showing that if $T \leq t \leq T + 2\pi$ satisfies
$$ \sum_{\gamma}\frac{1}{|t-\gamma|^{2}} \leq C\log^{2}T , $$
where $\gamma$ denotes the ordinates of non-trivial zeta zeros, then
$$  \log|\zeta(1/2+it)| = \Re \sum_{n \leq T} \frac{(\Lambda(n) / \log n)}{n^{1/2 + it}} \frac{\log(T/n)}{\log T} + O(1) , $$
where $\Lambda(n)$ denotes the von Mangoldt function. Afterwards we will show that the sum over $\gamma$ is indeed suitably small for most $T \leq t \leq T + 2\pi$, and that $\Re \sum_{n \leq T} \frac{(\Lambda(n) / \log n)}{n^{1/2 + it}} \frac{\log(T/n)}{\log T}$ can be replaced by $\Re \sum_{p \leq T} \frac{\log(T/p)/\log T}{p^{1/2+ it}}$ for most such $t$.

Since we assume that $\sum_{\gamma} 1/|t-\gamma|^{2} \leq C\log^{2}T$, and we are always assuming that $T + 2\pi \geq t \geq T$ is large, we may assume that the zeta function has no zeros or poles on the horizontal line extending from $1/2+it$ to (positive) infinity. Thus we have, assuming the Riemann Hypothesis, the following estimates:
\begin{eqnarray}
\log|\zeta(1/2+it)| & = & \Re\left(\sum_{n \leq T} \frac{(\Lambda(n) / \log n)}{n^{1/2 + it}} \frac{\log(T/n)}{\log T} - \frac{1}{\log T} \frac{\zeta'}{\zeta}(1/2+it) + \right. \nonumber \\
&& \left. + \frac{1}{\log T} \sum_{\gamma} \int_{1/2}^{\infty} \frac{T^{1/2+i\gamma-\sigma-it}}{(1/2+i\gamma-\sigma-it)^{2}} d\sigma + O\left(\frac{1}{\log T} \right) \right) ; \nonumber
\end{eqnarray}
$$ \Re \frac{\zeta'}{\zeta}(1/2+it) = -(1/2)\log T + O(1) . $$
These are essentially obtained on pages 4 and 5 of Soundararajan's article~\cite{soundmoments}, the first by integrating an explicit formula for $\zeta'/\zeta$ from $1/2+it$ to positive infinity along a horizontal line, and the second by taking real parts in the Hadamard (partial fraction) formula for $\zeta'/\zeta$. Inserting the second estimate in the first, we deduce that
$$ \log|\zeta(1/2+it)| = \Re \sum_{n \leq T} \frac{(\Lambda(n) / \log n)}{n^{1/2 + it}} \frac{\log(T/n)}{\log T} + O(1) + O\left(\frac{1}{\log T} (\sum_{\gamma} \frac{1}{|\gamma - t|^{2}}) \int_{1/2}^{\infty} T^{1/2-\sigma} d\sigma \right) . $$
The integral here is equal to $1/\log T$, so if we have $\sum_{\gamma} 1/|t-\gamma|^{2} \leq C\log^{2}T$ then the third term is $O(1)$, as claimed.

Next we define two ``good'' sets $\mathcal{G}^{(2)} \subseteq \mathcal{G}^{(1)} \subseteq [T,T+2\pi]$, by
$$ \mathcal{G}^{(1)} := \{T \leq t \leq T+2\pi : |t-\gamma| > 1/(K\log T) \; \forall \gamma \}, $$
$$ \mathcal{G}^{(2)} := \{t \in \mathcal{G}^{(1)} : \sum_{T-1 \leq \gamma \leq T +7} 1/|t-\gamma|^{2} \leq (C/2)\log^{2}T \}, $$
where $\gamma$ continues to run over ordinates of zeta zeros. We see
\begin{eqnarray}
\int_{\mathcal{G}^{(1)}} \sum_{T-1 \leq \gamma \leq T +7}\frac{1}{|t-\gamma|^{2}} dt = \sum_{T-1 \leq \gamma \leq T+7} \int_{\mathcal{G}^{(1)}} \frac{1}{|t-\gamma|^{2}} dt & \leq & 2 \sum_{T-1 \leq \gamma \leq T+7} \int_{1/(K\log T)}^{\infty} \frac{1}{t^{2}} dt \nonumber \\
& = & 2K\log T \sum_{T-1 \leq \gamma \leq T+7} 1 , \nonumber
\end{eqnarray}
and using standard estimates for the number of zeta zeros in a horizontal strip (as in e.g. Theorem 10.13 of Montgomery and Vaughan~\cite{mv}) we find this is $\ll K\log^{2}T$. Thus, provided $C$ was set sufficiently large in terms of $K$, the measure of $\mathcal{G}^{(1)} \backslash \mathcal{G}^{(2)}$ is at most $0.001\pi$. Similarly, provided $K$ was set sufficiently large the measure of $[T,T+2\pi] \backslash \mathcal{G}^{(1)}$ is at most $0.001\pi$ (since there are $\ll \log T$ zeros in that interval). If $t \in \mathcal{G}^{(2)}$ then $\sum_{\gamma} 1/|t-\gamma|^{2} \leq \sum_{T-1 \leq \gamma \leq T+7} 1/|t-\gamma|^{2} + \sum_{\gamma} 10/(1+|t-\gamma|^{2}) \leq C\log^{2}T$, and so we have this desired bound for all $T \leq t \leq T+2\pi$, except for a ``bad'' set of measure at most $0.002\pi$.

Finally we note that, since $\sum_{p \leq \sqrt{T}} (\log p)/p \ll \log T$, we have
\begin{eqnarray}
\sum_{n \leq T} \frac{(\Lambda(n) / \log n)}{n^{1/2 + it}} \frac{\log(T/n)}{\log T} & = & \sum_{p \leq T} \frac{1}{p^{1/2 + it}} \frac{\log(T/p)}{\log T} + \sum_{p \leq \sqrt{T}} \frac{1/2}{p^{1 + 2it}} \frac{\log(T/p^{2})}{\log T} + O(1) \nonumber \\
& = & \sum_{p \leq T} \frac{1}{p^{1/2 + it}} \frac{\log(T/p)}{\log T} + \sum_{p \leq \sqrt{T}} \frac{1/2}{p^{1 + 2it}} + O(1) ; \nonumber
\end{eqnarray}
so it will complete the proof of Proposition 1 if we can show that
$$ \text{meas}\{T \leq t \leq T+2\pi : \left| \sum_{p \leq \sqrt{T}} \frac{1}{p^{1 + 2it}}  \right| > C\} \leq 0.008\pi . $$
But this is an easy consequence of the fact that $ \int_{0}^{2\pi} \left| \sum_{p \leq \sqrt{T}} \frac{p^{-2iT}}{p^{1 +2it}} \right|^{2} dt \ll 1 $,
which follows from a suitable form of Plancherel's identity (as in e.g. equation (5.26) of Montgomery and Vaughan~\cite{mv}).
\begin{flushright}
Q.E.D.
\end{flushright}

\section{The probabilistic part : Proof of Proposition 2}

\subsection{The upper bound}
For the sake of concision, let us write
$$ V(p,h) := \frac{U_{p} p^{-ih}}{p^{1/\log T}} \frac{\log(T/p)}{\log T} + \textbf{1}_{p \leq \sqrt{T}} \frac{1}{2} \frac{U_{p}^{2} p^{-2ih}}{p^{1/2+2/\log T}} \frac{\log(T/p^{2})}{\log T} , \;\;\;\;\; p \leq T, \;\;\; 0 \leq h \leq 2\pi , $$
where $\textbf{1}$ denotes the indicator function. To prove the upper bound in Proposition 2, we need to show that there exists an absolute constant $C > 0$ such that
$$ \p(\max_{0 \leq h \leq 2\pi} \Re \sum_{p \leq T} \frac{V(p,h)}{p^{1/2}} > \log\log T - (1/4)\log\log\log T + C\sqrt{\log\log\log T}) = o(1) $$
as $T \rightarrow \infty$.

Using the union bound, this probability is certainly at most
$$ \sum_{0 \leq j \leq 2\pi\log T} \p\left(\max_{\frac{j}{\log T} \leq h \leq \frac{j+1}{\log T}} \Re \sum_{p \leq T} \frac{V(p,h)}{p^{1/2}} > \log\log T - (1/4)\log\log\log T + C\sqrt{\log\log\log T}\right) , $$
and since $U_{p} p^{-ij/\log T}$ has the same distribution as $U_{p}$, for any $j$, this is
$$ \ll \log T \cdot \p\left(\max_{0 \leq h \leq \frac{1}{\log T}} \Re \sum_{p \leq T} \frac{V(p,h)}{p^{1/2}} > \log\log T - (1/4)\log\log\log T + C\sqrt{\log\log\log T}\right) . $$
So it will suffice to show that the probability above is $o(1/\log T)$ as $T \rightarrow \infty$.

The random sums $\Re \sum_{p \leq T} \frac{V(p,h_1)}{p^{1/2}}, \Re \sum_{p \leq T} \frac{V(p,h_2)}{p^{1/2}}$ are almost perfectly correlated when $|h_1 - h_2| \leq 1/\log T$, so one might expect the probability above to be roughly the same as $ \p\left(\Re \sum_{p \leq T} \frac{V(p,0)}{p^{1/2}} > \log\log T - (1/4)\log\log\log T + C\sqrt{\log\log\log T}\right)$. This is essentially true, and in fact we will shortly prove the following result:
\begin{lem1}
If $T$ is large, and $C > 0$ is a sufficiently large absolute constant, then
\begin{eqnarray}
&& \p\left(\max_{0 \leq h \leq \frac{1}{\log T}} \Re \sum_{p \leq T} \frac{V(p,h)}{p^{1/2}} > \log\log T - (1/4)\log\log\log T + C\sqrt{\log\log\log T}\right) \nonumber \\
& \leq & \p\left(\Re \sum_{p \leq T} \frac{V(p,0)}{p^{1/2}} > \log\log T - (1/4)\log\log\log T + \sqrt{\log\log\log T}\right) + O(\frac{1}{\log T \log\log T}) . \nonumber
\end{eqnarray}
\end{lem1}

Now we need a bound for the probability that $\Re \sum_{p \leq T} V(p,0)/p^{1/2}$, which is simply a sum of independent random variables, is large. There are many such bounds available, but the next result, which follows from Theorem 3.3 of Talagrand~\cite{talagrand} (by choosing $u=t/\sigma^{2}$ to bound the infimum of characteristic functions there), is slightly sharper than most of those and will allow us\footnote{Recall that if $Z$ is a standard normal random variable, and $z$ is large, then $\p(Z > z) \ll (1/z) e^{-z^{2}/2}$. Most tail probability bounds for sums of independent random variables have an exponential component like this, but do not include the multiplier $1/z$. Talagrand's result supplies such a multiplier, which will let us obtain a probability bound $o(1/\log T)$ at a slightly lower threshold.} to have the subtracted term $-(1/4)\log\log\log T$.
\begin{tail}[Talagrand, 1995]
There exists an absolute constant $K > 0$ such that the following is true. Suppose $X_{1},...,X_{n}$ are independent, mean zero random variables, and suppose $B > 0$ is such that $|X_{i}| \leq B$ almost surely, for all $i$. Set $\sigma^{2} := \E(\sum_{1 \leq i \leq n} X_{i})^{2} = \sum_{1 \leq i \leq n} \E X_{i}^{2}$. Then for any $0 \leq t \leq \sigma^{2}/KB$ we have
$$ \p(\sum_{1 \leq i \leq n} X_{i} \geq t) \ll \frac{1}{1+(t/\sigma)} \exp\left(-t^{2}/(2\sigma^{2}) + O(|t/\sigma^{2}|^{3} \sum_{1 \leq i \leq n} \E|X_{i}|^{3})\right) . $$
\end{tail}

Applying Tail Inequality 1 to $\sum_{1000K^{2} \leq p \leq T} \Re(V(p,0)/p^{1/2})$, we see that we can take $B = 1/(10K)$ (say), and we have
\begin{eqnarray}
\sigma^{2} = \sum_{1000K^{2} \leq p \leq T} \E (\Re \frac{V(p,0)}{p^{1/2}})^{2} & = & \sum_{1000K^{2} \leq p \leq T} \frac{1}{p} (\E (\Re U_{p})^{2} + O(\frac{\log p}{\log T}) + O(\frac{1}{\sqrt{p}})) \nonumber \\
& = & (1/2) \sum_{1000K^{2} \leq p \leq T} \frac{1}{p} + O(1) = (1/2)\log\log T + O(1) \nonumber
\end{eqnarray}
and $\sum_{1000K^{2} \leq p \leq T} \E|\Re(V(p,0)/p^{1/2})|^{3} \ll \sum_{1000K^{2} \leq p \leq T} 1/p^{3/2} \ll 1$. See the first appendix, below, for some similar variance calculations. Thus we have
\begin{eqnarray}
&& \p\left(\Re \sum_{p \leq T} \frac{V(p,0)}{p^{1/2}} > \log\log T - (1/4)\log\log\log T + \sqrt{\log\log\log T}\right) \nonumber \\
& \leq & \p\left(\Re(\sum_{1000K^{2} \leq p \leq T} \frac{V(p,0)}{p^{1/2}}) > \log\log T - (1/4)\log\log\log T + (1/2)\sqrt{\log\log\log T}\right) \nonumber \\
& \ll & \frac{1}{\sqrt{\log\log T}} \exp\left(\frac{-(\log\log T - (1/4)\log\log\log T + (1/2)\sqrt{\log\log\log T})^{2}}{\log\log T + O(1)} + O(1)\right) \nonumber \\
& = & \frac{1}{\sqrt{\log\log T}} \exp\left(-\log\log T + (1/2)\log\log\log T - \sqrt{\log\log\log T} + O(1) \right) , \nonumber
\end{eqnarray}
provided that $T$ is large enough. This bound is $o(1/\log T)$, which suffices to establish the upper bound in Proposition 2 provided we can prove Lemma 1. 

\vspace{12pt}
To prove the lemma we will use a {\em chaining argument}, whereby we approximate $\max_{0 \leq h \leq \frac{1}{\log T}} \Re \sum_{p \leq T} \frac{V(p,h)}{p^{1/2}}$ by the maximum over increasingly sparse discrete sets of points $h$. This argument would be easier if the random sums $\Re \sum_{p \leq T} \frac{V(p,h)}{p^{1/2}}$ were Gaussian random variables, and in $\S\S 3.2-3.3$, where we prove lower bounds, we will have to use a central limit theorem to pass to that situation. However, here we will make do without that step.

Let $\mathcal{H}_{0} := \{0\}$, and for $1 \leq k \leq \log T$ let $\mathcal{H}_{k} := \{ i/(2^{k}\log T) : 0 \leq i \leq 2^{k}\}$ (so that $\mathcal{H}_{k-1} \subset \mathcal{H}_{k} \; \forall k$), and note that we certainly have
\begin{eqnarray}
\max_{0 \leq h \leq \frac{1}{\log T}} \Re \sum_{p \leq T} \frac{V(p,h)}{p^{1/2}} & = & \max_{h \in \mathcal{H}_{[\log T]}} \Re \sum_{p \leq T} \frac{V(p,h)}{p^{1/2}} + O(1) \nonumber \\
& \leq & \max_{h \in \mathcal{H}_{[\log T]}} \Re(\sum_{p \leq T^{1/\log\log T}} \frac{V(p,h)}{p^{1/2}}) + \max_{h \in \mathcal{H}_{[\log T]}} \Re(\sum_{\substack{T^{1/\log\log T} \\ < p \leq T}} \frac{V(p,h)}{p^{1/2}}) + O(1) . \nonumber
\end{eqnarray}
We split into two pieces here because, as the reader will soon see, the part of the sum over small primes is more highly correlated at short distances (and also contributes most of the expected size of the maximum), and so it is easier to handle. Now given a point $h \in \mathcal{H}_{[\log T]}$, define $h^{(k)} := \max\{t \leq h : t \in \mathcal{H}_{k}\}$, and note that $\Re \sum_{p \leq T^{1/\log\log T}} \frac{V(p,h)}{p^{1/2}}$ is
\begin{eqnarray}
& = & \Re(\sum_{p \leq T^{1/\log\log T}} \frac{V(p,0)}{p^{1/2}}) + \sum_{1 \leq k \leq \log T} \left(\Re(\sum_{p \leq T^{1/\log\log T}} \frac{V(p,h^{(k)})}{p^{1/2}}) - \Re(\sum_{p \leq T^{1/\log\log T}} \frac{V(p,h^{(k-1)})}{p^{1/2}}) \right) \nonumber \\
& \leq & \Re(\sum_{p \leq T^{1/\log\log T}} \frac{V(p,0)}{p^{1/2}}) + \sum_{1 \leq k \leq \log T} \max_{h \in \mathcal{H}_{k}} \left(\Re(\sum_{p \leq T^{1/\log\log T}} \frac{V(p,h) - V(p,h^{(k-1)})}{p^{1/2}}) \right) . \nonumber
\end{eqnarray}
This is the chain decomposition that will drive our argument. Indeed, unless there is some $k$ such that $\max_{h \in \mathcal{H}_{k}} \Re(\sum_{p \leq T^{1/\log\log T}} \frac{V(p,h) - V(p,h^{(k-1)})}{p^{1/2}}) > k^{0.9} 2^{-k}$ (say) the contribution from $\sum_{1 \leq k \leq \log T}$ will clearly be $O(1)$, and therefore for any $u \in \R$ we have
\begin{eqnarray}\label{chainbound}
&& \p\left(\max_{0 \leq h \leq \frac{1}{\log T}} \Re \sum_{p \leq T} \frac{V(p,h)}{p^{1/2}} > u \right) \nonumber \\
& \leq &  \p\left(\Re(\sum_{p \leq T^{1/\log\log T}} \frac{V(p,0)}{p^{1/2}}) + \max_{h \in \mathcal{H}_{[\log T]}} \Re(\sum_{\substack{T^{1/\log\log T} \\ < p \leq T}} \frac{V(p,h)}{p^{1/2}})  > u - O(1) \right) + \nonumber \\
&& + \sum_{1 \leq k \leq \log T} \sum_{h \in \mathcal{H}_{k}} \p\left(\Re(\sum_{p \leq T^{1/\log\log T}} \frac{V(p,h) - V(p,h^{(k-1)})}{p^{1/2}}) > k^{0.9} 2^{-k} \right) .
\end{eqnarray}

Next, for any $h \in \mathcal{H}_{k-1}$ the sum $\Re(\sum_{p \leq T^{1/\log\log T}} \frac{V(p,h) - V(p,h^{(k-1)})}{p^{1/2}})$ vanishes, and for any $h \in \mathcal{H}_{k} \backslash \mathcal{H}_{k-1}$ it has mean zero and variance
\begin{eqnarray}
\sum_{p \leq T^{1/\log\log T}} \E(\Re\frac{V(p,h) - V(p,h^{(k-1)})}{p^{1/2}})^{2} & \asymp & \sum_{p \leq T^{1/\log\log T}} \frac{1}{p} \E(\Re U_{p} (p^{-ih}-p^{-ih^{(k-1)}}))^{2} \nonumber \\
& \asymp & \sum_{p \leq T^{1/\log\log T}} \frac{1}{p} |1-p^{i(h-h^{(k-1)})}|^{2} \nonumber \\
& \asymp & \sum_{p \leq T^{1/\log\log T}} \frac{1}{p} \frac{\log^{2}p}{2^{2k} \log^{2}T} \asymp \frac{1}{(2^{k} \log\log T)^{2}} , \nonumber
\end{eqnarray}
since $|1-p^{i(h-h^{(k-1)})}| \asymp |1-p^{i/(2^{k}\log T)}|$ for all $h \in \mathcal{H}_{k} \backslash \mathcal{H}_{k-1}$, and $\sum_{p \leq x} \frac{\log^{2}p}{p} \asymp \log^{2}x$. Using Tail Inequality 1 (with $B=100/(2^{k}\log T)$, say), it follows that for all $h \in \mathcal{H}_{k}$,
$$ \p\left(\Re(\sum_{p \leq T^{1/\log\log T}} \frac{V(p,h) - V(p,h^{(k-1)})}{p^{1/2}}) > k^{0.9} 2^{-k} \right) \ll \exp(- c k^{1.8} (\log\log T)^{2}) \ll \frac{2^{-k}}{\log^{3}T} , $$
where $c > 0$ is a small absolute constant. Thus the sums in (\ref{chainbound}) are collectively of size $O(1/\log^{2}T)$, which is more than acceptable for Lemma 1. A similar chaining argument can be applied to $\max_{h \in \mathcal{H}_{[\log T]}} \Re(\sum_{T^{1/\log\log T} < p \leq T} \frac{V(p,h)}{p^{1/2}})$, but in this case the variance of the pieces is $O(2^{-2k})$ rather than $O(2^{-2k} (\log\log T)^{-2})$, and so we can only immediately handle large values of $k$ (say $\log\log\log T \leq k \leq \log T$). We can deduce that, for any $u \in \R$, the probability $\p\left(\max_{0 \leq h \leq \frac{1}{\log T}} \Re \sum_{p \leq T} \frac{V(p,h)}{p^{1/2}} > u \right)$ is
$$ \leq \p\left(\Re(\sum_{p \leq T^{1/\log\log T}} \frac{V(p,0)}{p^{1/2}}) + \max_{h \in \mathcal{H}_{[\log\log\log T]}} \Re(\sum_{\substack{T^{1/\log\log T} \\ < p \leq T}} \frac{V(p,h)}{p^{1/2}})  > u - O(1) \right) + O(\frac{1}{\log^{2}T}) . $$

Finally set $u = \log\log T - (1/4)\log\log\log T + C\sqrt{\log\log\log T}$, as in Lemma 1, and notice if $\Re(\sum_{p \leq T^{1/\log\log T}} \frac{V(p,0)}{p^{1/2}}) + \max_{h \in \mathcal{H}_{[\log\log\log T]}} \Re(\sum_{T^{1/\log\log T} < p \leq T} \frac{V(p,h)}{p^{1/2}})  > u - O(1)$ then one of the following must occur:
\begin{itemize}
\item $\Re(\sum_{p \leq T^{1/\log\log T}} \frac{V(p,0)}{p^{1/2}}) + \max_{h \in \mathcal{H}_{[\log\log\log T]}} \Re(\sum_{T^{1/\log\log T} < p \leq T} \frac{V(p,h)}{p^{1/2}})  > u - O(1)$, and also
$$ \max_{h \in \mathcal{H}_{[\log\log\log T]}} \Re(\sum_{T^{1/\log\log T} < p \leq T} \frac{V(p,h)}{p^{1/2}}) > (C/2)\log\log\log T ; $$

\item $\Re(\sum_{p \leq T^{1/\log\log T}} \frac{V(p,0)}{p^{1/2}}) > u - O(1) - (C/2)\log\log\log T$, and also
$$ \max_{h \in \mathcal{H}_{[\log\log\log T]}} \left( \Re(\sum_{T^{1/\log\log T} < p \leq T} \frac{V(p,h)}{p^{1/2}}) - \Re(\sum_{T^{1/\log\log T} < p \leq T} \frac{V(p,0)}{p^{1/2}}) \right)  > (C/2)\sqrt{\log\log\log T} ; $$

\item $\Re(\sum_{p \leq T^{1/\log\log T}} \frac{V(p,0)}{p^{1/2}}) > u - O(1) - (C/2)\log\log\log T$, and also
$$ \Re(\sum_{p \leq T} \frac{V(p,0)}{p^{1/2}})  > u - (C/2)\sqrt{\log\log\log T} - O(1) . $$
\end{itemize}
The probability of the third event here is at most as large as the probability in the conclusion of Lemma 1, so to prove the lemma it will suffice to show that each of the other events has probability $O(1/(\log T \log\log T))$. This follows using the independence of the sums over $p \leq T^{1/\log\log T}$ and over $T^{1/\log\log T} < p \leq T$, and using the union bound. For example, the probability of the second event is
\begin{eqnarray}
& \leq & \p\left( \Re(\sum_{p \leq T^{1/\log\log T}} \frac{V(p,0)}{p^{1/2}}) > u - O(1) - (C/2)\log\log\log T \right) \cdot \nonumber \\
&& \cdot \sum_{h \in \mathcal{H}_{[\log\log\log T]}} \p\left( \Re(\sum_{T^{1/\log\log T} < p \leq T} \frac{V(p,h) - V(p,0)}{p^{1/2}}) > (C/2)\sqrt{\log\log\log T} \right) \nonumber \\
& \ll & \exp(\frac{-(u-O(1)-(C/2)\log\log\log T)^{2}}{\log\log T}) \cdot \log\log T \cdot \exp(-c((C/2)\sqrt{\log\log\log T})^{2}) \nonumber \\
& \ll & \log\log T \cdot \exp(-\log\log T + 2C\log\log\log T -c(C/2)^{2}\log\log\log T) , \nonumber
\end{eqnarray}
say, where $c > 0$ is again a small absolute constant, and we used the facts that $\#\mathcal{H}_{[\log\log\log T]} \ll \log\log T$ and $\E (\Re \sum_{T^{1/\log\log T} < p \leq T} \frac{V(p,h) - V(p,0)}{p^{1/2}} )^{2} = O(1)$ (and Tail Inequality 1, with $B=100/T^{1/2\log\log T}$). This is indeed $O(1/(\log T \log\log T))$, provided that $C > 0$ was chosen sufficiently large, and a similar argument (using the fact that $\E (\Re \sum_{T^{1/\log\log T} < p \leq T} \frac{V(p,h)}{p^{1/2}} )^{2} \sim (1/2)\log\log\log T$) applies to the probability of the first event.
\begin{flushright}
Q.E.D.
\end{flushright}

Apart from the chaining arguments used to prove Lemma 1, the proof of the upper bound in Proposition 2 was essentially just an application of the union bound on the scale of $1/\log T$. If the random sums $\Re \sum_{p \leq T} \frac{V(p,h_1)}{p^{1/2}}, \Re \sum_{p \leq T} \frac{V(p,h_2)}{p^{1/2}}$ behaved independently when $|h_1 - h_2| \gg 1/\log T$ then one would expect such an argument to be quite sharp, since if $X_1,...,X_n$ are independent then
$$ \p(\max_{1 \leq i \leq n} X_{i} > u) = 1 - \p(X_{i} \leq u \; \forall i) = 1 - \prod_{i=1}^{n}(1- \p(X_{i} > u)) \approx \min\{1, \sum_{i=1}^{n} \p(X_{i} > u)\} , $$
which follows since $1- \p(X_{i} > u) \approx \exp(-\p(X_{i} > u))$. This is why the upper bound in Proposition 2, with the subtracted term $-(1/4)\log\log\log T$, would be sharp in the ``independent at distance $1/\log T$'' case, but since the sums $\Re \sum_{p \leq T} \frac{V(p,h_1)}{p^{1/2}}$ are actually logarithmically correlated we expect their maximum to be a little smaller (since we have ``fewer independent tries at obtaining a large value'').

We also note some remarks made by Fyodorov and Keating about the long range maximum $ \max_{0 \leq t \leq T} |\zeta(1/2+it)|$ studied by Farmer, Gonek and Hughes~\cite{farmergonekhughes}. At the end of $\S 2.5$ of their paper~\cite{fyodkeat}, Fyodorov and Keating observe: ``the tail of the distribution [of $- \log\max_{T \leq t \leq T + 2\pi} |\zeta(1/2+it)| + \log\log T - (3/4)\log\log\log T$] that we predict for much shorter ranges decays like $|x|e^{x}$ as $x \rightarrow -\infty$; that is, the exponential is linear rather than quadratic. If this were to persist... it would suggest that $\zeta(1/2+it)$ may take much larger values than... the Farmer-Gonek-Hughes conjecture''. They add that ``there are several reasons for thinking this unlikely''. The foregoing calculations easily imply that, in our model, the {\em long} range tail decays like a {\em quadratic} exponential, and similar considerations quite possibly apply to Fyodorov and Keating's random matrix model.

\subsection{Preliminary calculations for the lower bound}
In this subsection we make some preliminary modifications to the collection of random variables
$$ \Re \sum_{p \leq T} \frac{U_{p}}{p^{1/2+ih}} \frac{\log(T/p)}{\log T} , \;\;\;\;\; h \in \mathcal{H} $$
in the lower bound part of Proposition 2. At the end of the subsection we will give an overview of the reasons for making these modifications.

Firstly let $1 \ll E=E(T) \ll \sqrt{\log T}$ be a parameter, whose value will be fixed in $\S 3.3$. (In fact we will end up taking $E = \sqrt{\log\log T} (\log\log\log T)^{2}$). We claim that there exists some $0 \leq z \leq E/\log T$ such that at least $1.98 \pi (\log T)/E$ of the points
$$ z + \frac{iE}{\log T}, \;\;\;\;\; 0 \leq i \leq (2\pi\log T)/E - 1 $$
belong to the ``good'' set $\mathcal{H}$. Indeed, this follows immediately when we observe that, if $\textbf{1}$ denotes the indicator function,
$$ \int_{0}^{E/\log T} \sum_{0 \leq i \leq (2\pi\log T)/E - 1} \textbf{1}_{z + iE/\log T \in \mathcal{H}} dz \geq \text{meas}(\mathcal{H}) - \frac{E}{\log T} \geq 1.98 \pi . $$
We choose such a value of $z$ and let $\mathcal{H}^{*} := \mathcal{H} \cap \{z + iE/\log T : 0 \leq i \leq (2\pi\log T)/E - 1\}$, a discretisation of the set $\mathcal{H}$ that will be more convenient to work with.

Now let $2 \leq y=y(T) \ll e^{(\log\log T)^{1000}}$ be a further parameter, whose value will also be fixed in $\S 3.3$. (In fact we will end up taking $y = e^{(\log\log T)^{2} (\log\log\log T)^{2}}$). As we will explain shortly, for technical reasons (see below, and also the calculations and discussion following Lemma 3 in $\S 3.3$) we need to remove the primes smaller than $y$ from our random sums. To account for the error that arises in doing this we can take a very crude approach: for any fixed $h \in \mathcal{H}^{*}$ we have
$$ \E \left(\Re \sum_{p < y} \frac{U_{p}}{p^{1/2+ih}} \frac{\log(T/p)}{\log T} \right)^{2} \ll \log\log y \ll \log\log\log T , $$
using the variance calculations in our first appendix (with $P=2$ and $Q=y$), and therefore by Chebychev's inequality we have
$$ \p(\left|\Re \sum_{p < y} \frac{U_{p}}{p^{1/2+ih}} \frac{\log(T/p)}{\log T} \right| \geq (\log\log\log T)^{3/4}) \ll (\log\log\log T)^{-1/2} = o(1) . $$
Since the random variables $(U_{p})_{p < y}$ are independent of $(U_{p})_{y \leq p \leq T}$, we see that to prove the lower bound in Proposition 2 it will suffice to show that
$$ \p(\max_{h \in \mathcal{H}^{*}} \Re \sum_{y \leq p \leq T} \frac{U_{p}}{p^{1/2+ih}} \frac{\log(T/p)}{\log T} \geq \log\log T - 2\log\log\log T - (\log\log\log T)^{3/4}) = 1 - o(1) . $$

Next we let
$$ Y(h) := \frac{\Re \sum_{y \leq p \leq T} \frac{U_{p}}{p^{1/2+ih}} \frac{\log(T/p)}{\log T} }{\sqrt{\frac{1}{2} \sum_{y \leq p \leq T} \frac{1}{p} \frac{\log^{2}(T/p)}{\log^{2}T}}} , \;\;\;\;\; h \in \mathcal{H}^{*} . $$
Using the variance calculations in our first appendix (with $P=y$ and $Q=T$), we see the $Y(h)$ are mean zero, variance one, real-valued random variables. Moreover, if $h_{1} \neq h_{2} \in \mathcal{H}^{*}$ (so that, in particular, $1/\log T \leq |h_{1}-h_{2}| \leq 2\pi$) then the calculations in our first appendix supply the following covariance estimate:
\begin{equation}\label{correst}
\E Y(h_1)Y(h_2) = \left\{ \begin{array}{ll}
     1 - \frac{\log(|h_{1}-h_{2}|\log T) + O(1)}{\log\log T - \log\log y + O(1)} & \text{if} \; 1/\log T \leq |h_{1}-h_{2}| \leq 1/\log y   \\
     O(\frac{1}{|h_{1}-h_{2}| \log y \log\log T}) & \text{otherwise.}
\end{array} \right.
\end{equation}
Here we used the fact, justified in our first appendix, that
$$  \sum_{y \leq p \leq T} \frac{1}{p} \frac{\log^{2}(T/p)}{\log^{2}T} = \log\log T - \log\log y + O(1) \gg \log\log T , $$
bearing in mind that $y \ll e^{(\log\log T)^{1000}}$. We also remark that the random variables $Y(h)$ are {\em stationary}, in other words the covariance $\E Y(h_1)Y(h_2)$ only depends on the distance $|h_{1} - h_{2}|$. This is visibly true at the level of the estimates (\ref{correst}), if one ignores ``big Oh'' terms, and in fact it is exactly true, since we actually have
$$ \E Y(h_1)Y(h_2) = \frac{\frac{1}{2} \sum_{y \leq p \leq T} \frac{\cos((h_{1}-h_{2})\log p)}{p} \frac{\log^{2}(T/p)}{\log^{2}T}}{\frac{1}{2} \sum_{y \leq p \leq T} \frac{1}{p} \frac{\log^{2}(T/p)}{\log^{2}T}} ,  $$
as shown in our first appendix. The exact stationarity of the $Y(h)$ isn't really necessary for the analysis in $\S 3.3$, (see the author's paper~\cite{harpergp}, where some very similar random variables without this property are treated), but it is quite convenient.

Finally, using a multivariate central limit theorem (as explained in our second appendix) we can replace the random variables $Y(h)$, $h \in \mathcal{H}^{*}$ by Gaussian random variables with the same means and covariance matrix, provided that $y \geq \log^{7}T$, say (so that none of the summands in the definition of $Y(h)$ is too large relative to $\#\mathcal{H}^{*}$). We summarise the state of affairs we have reached in the following lemma.
\begin{lem2}
Suppose $1 \ll E \ll \sqrt{\log T}$ and $\log^{7}T \ll y \ll e^{(\log\log T)^{1000}}$, and let $\mathcal{H}^{*}$ be as above. Let $Z(h)$, $h \in \mathcal{H}^{*}$ be a collection of mean zero, variance one, jointly normal random variables with the same covariances as the random variables $Y(h)$ described above. Then to prove the lower bound in Proposition 2, it will suffice to show that
$$ \p(\max_{h \in \mathcal{H}^{*}} Z(h) \geq \frac{\log\log T - 2\log\log\log T - (\log\log\log T)^{3/4} + 1}{\sqrt{(1/2)(\log\log T - \log\log y + O(1))}}) = 1 - o(1) . $$
\end{lem2}

Our reason for transitioning to Gaussian random variables is because all information about their dependencies is contained in their covariances, and many tools exist for analysing their behaviour, as we shall see in $\S 3.3$. We introduced the parameter $y$ to make the random variables $Z(h)$ less correlated at distant values of $h$ (note the factor $\log y$ in the denominator of the second estimate in (\ref{correst})). This will be necessary to complete the proof of Proposition 2 as stated, although one could obtain a weaker result without doing this, at least for the Gaussian $Z(h)$. The parameter $E$ spaces out the points in $\mathcal{H}^{*}$ so that fewer of them are very close together, and therefore fewer of the $Z(h)$ are very highly correlated. As the reader will see in $\S 3.3$, this seems to be quite essential to obtain a successful conclusion.

\subsection{The lower bound}
Our main tool for establishing the lower bound in Proposition 2 is the following, which is a slight adaptation\footnote{Theorem 1 of \cite{harpergp} deals with the case where $\mathcal{N} = \{1,2,...,n\}$. However, to treat the more general case one can simply replace the sum over $1 \leq m \leq n$ in Proposition 1 of \cite{harpergp} by the corresponding sum over $m \in \mathcal{N}$.} of Theorem 1 from the author's paper~\cite{harpergp}.
\begin{lowerb}[Harper, 2013]
Let $\{Z(t_{i})\}_{1 \leq i \leq n}$ be jointly multivariate normal random variables, each with mean zero and variance 1. Suppose that the sequence is {\em stationary}, i.e. that $\E Z(t_{i})Z(t_{j}) = r(|i-j|)$ for some function $r$. Finally let $u \geq 1$, and suppose that:
\begin{itemize}
\item $r(m)$ is a decreasing non-negative function;

\item $r(1)(1+2u^{-2})$ is at most $1$.
\end{itemize}
Then for any subset $\mathcal{N} \subseteq \{1,2,...,n\}$, the probability $\p(\max_{i \in \mathcal{N}} Z(t_{i}) > u)$ is
$$ \geq (\#\mathcal{N}) \frac{e^{-u^{2}/2}}{40 u} \min\left\{1,\sqrt{\frac{1-r(1)}{u^{2}r(1)}}\right\} \prod_{j=1}^{n-1} \Phi\left(u\sqrt{1-r(j)} \left(1+O\left(\frac{1}{u^{2}(1-r(j))}\right) \right) \right), $$
where $\Phi(z) := (1/\sqrt{2\pi}) \int_{-\infty}^{z} e^{-t^{2}/2} dt$ is the standard normal distribution function.
\end{lowerb}

\vspace{12pt}
To get an idea of what Lower Bound 1 says, we will first apply it to the random variables $Z(h)$ in Lemma 2 without carefully checking the conditions, and ignoring the ``big Oh'' error terms in the covariances (\ref{correst}) and in the statement of the theorem. Afterwards we will explain how to apply the theorem properly.

Since $\mathcal{H}^{*}$ is a set of points of the form $z + (iE)/\log T$, we can clearly reparametrise the $(Z(h))_{h \in \mathcal{H}^{*}}$ by integers $i$, as in Lower Bound 1. Moreover, as remarked in our Preliminary Calculations the covariances $\E Z(h_1)Z(h_2)$ only depend on $|h_1 - h_2|$, so these random variables are indeed stationary. Thus if $\mathcal{I}$ is some set of indices $i$ that are contained in an interval of length $(\log T)/(E\log y)$, so that $|h_i - h_j| \leq 1/\log y$ for the corresponding points $h_{i} = z + (iE)/\log T, h_{j} = z + (jE)/\log T$, we have by (\ref{correst}) that
$$ 1 - r(j) := 1 - \E Z(h)Z(h + \frac{jE}{\log T}) \approx \frac{\log((jE/\log T)\log T)}{\log\log T - \log\log y} = \frac{\log(jE)}{\log\log T - \log\log y} , $$
and therefore $\p(\max_{i \in \mathcal{I}} Z(h_{i}) > \sqrt{2(\log\log T - \log\log y)})$ is
$$ \gtrsim (\#\mathcal{I}) \frac{\log y \sqrt{\log E}}{\log T (\log\log T)^{3/2}} \prod_{j=1}^{[(\log T)/(E\log y)]} \Phi(\sqrt{2\log(jE)}) . $$

To estimate the product, we note that if $z \geq 1$ then
\begin{eqnarray}
\Phi(z) = 1 - \frac{1}{\sqrt{2\pi}} \int_{z}^{\infty} e^{-t^{2}/2} dt \geq 1 - \frac{1}{z\sqrt{2\pi}} \int_{z}^{\infty} t e^{-t^{2}/2} dt & = & 1 - \frac{e^{-z^{2}/2}}{z\sqrt{2\pi}} \nonumber \\
& \geq & \exp(-e^{-z^{2}/2}/z), \nonumber
\end{eqnarray}
and therefore (since $E \gg 1$, so we certainly have $\sqrt{2\log(jE)} \geq 1$ for all $j$)
$$ \prod_{j=1}^{\left[\frac{\log T}{E\log y}\right]} \Phi(\sqrt{2\log(jE)}) \geq \exp\left(\sum_{j=1}^{\left[\frac{\log T}{E\log y}\right]} \frac{-1}{jE \sqrt{2\log(jE)}} \right) \geq \exp(- \frac{1}{E} \sum_{j=1}^{[\log T]} \frac{1}{j\sqrt{1 + \log j}} ) . $$
The sum here is $\ll \sqrt{\log\log T}$, so provided that $E \gg \sqrt{\log\log T}$ we will have
$$ \p(\max_{i \in \mathcal{I}} Z(h_{i}) > \sqrt{2(\log\log T - \log\log y)}) \gtrsim (\#\mathcal{I}) \frac{\log y \sqrt{\log E}}{\log T (\log\log T)^{3/2}} . $$
If the set $\mathcal{I}$ contains $\gg (\log T)/(E\log y)$ points then this lower bound will be ``not too small'', in a sense that will become clear shortly. Note that it is extremely important that the parameter $E$ be present here, and be large enough to compensate for the sum over $j$, since otherwise the lower bound becomes smaller by an exponential factor.

Actually it is not difficult to perform the foregoing calculations rigorously, provided we assume that $\mathcal{I}$ is contained in an interval of length $(\log T)/(K E\log y)$ for a suitable large constant $K > 0$. To make things rigorous we need to check that the function
$$ r(j) = \E Z(h)Z(h+\frac{jE}{\log T}) = \frac{\sum_{y \leq p \leq T} \frac{\cos((jE/\log T)\log p)}{p} \frac{\log^{2}(T/p)}{\log^{2}T}}{\sum_{y \leq p \leq T} \frac{1}{p} \frac{\log^{2}(T/p)}{\log^{2}T}} , \;\;\; 1 \leq j \leq (\log T)/(K E \log y) $$
is decreasing and non-negative; that $r(1)(1+1/(\log\log T - \log\log y)) \leq 1$; and that the ``big Oh'' error terms in the statement of Lower Bound 1, and in the correlations (\ref{correst}), do not alter the calculations. All of these things are straightforward to check using (\ref{correst}), except for the condition that $r(j)$ is decreasing and non-negative, which follows from the discussion at the end of our first appendix (with $P=y$ and $Q=T$) provided that $K$ is large enough and $\sqrt{\log y} \gg \log\log T$.

Again, we shall summarise the state of affairs we have reached as a lemma.
\begin{lem3}
Suppose that $\sqrt{\log\log T} \ll E \ll \sqrt{\log T}$ and $(\log\log T)^{2} \ll \log y \ll (\log\log T)^{1000}$, and let $Z(h)$, $h \in \mathcal{H}^{*}$ be the corresponding collection of mean zero, variance one, jointly multivariate normal random variables from Lemma 2.

Suppose that $\mathcal{I}$ is a set of integers contained in an interval of length $(\log T)/(KE \log y)$, such that $\#\mathcal{I} \geq (\log T)/(2KE\log y)$ (say) and such that
$$ h_{i} :=  z + \frac{iE}{\log T} \in \mathcal{H}^{*} \;\;\;\;\; \forall i \in \mathcal{I} . $$
Then $ \p(\max_{i \in \mathcal{I}} Z(h_{i}) > \sqrt{2(\log\log T - \log\log y)}) \gg \frac{\sqrt{\log E}}{E (\log\log T)^{3/2}} $, where the implicit constant is absolute.
\end{lem3}

\vspace{12pt}
Now recall, from our Preliminary Calculations, that $\mathcal{H}^{*}$ is a certain subset of $\{z + iE/\log T : 0 \leq i \leq (2\pi\log T)/E - 1\}$, and that $\mathcal{H}^{*}$ contains at least $1.98\pi (\log T)/E$ points. Thus if we let
$$ \mathcal{H}^{*}_{\text{even}} := \mathcal{H}^{*} \cap \bigcup_{\substack{0 \leq j \leq 2\pi K \log y , \\ j \; \text{even}}} \{z + iE/\log T : \frac{j \log T}{KE\log y} \leq i \leq \frac{(j+1) \log T}{KE\log y} \} , $$
then $\mathcal{H}^{*}_{\text{even}}$ must contain at least $0.95\pi (\log T)/E$ points, say, since the complementary set $\mathcal{H}^{*}_{\text{odd}}$ (where $j$ runs over odd integers) must satisfy
$$ \#\mathcal{H}^{*}_{\text{odd}} \leq (\pi K \log y + O(1))((\log T)/(K E\log y) + O(1)) \leq 1.03\pi (\log T)/E . $$
In particular, we may find sets $(\mathcal{I}_{k})_{1 \leq k \leq \log y}$, each a subset of $\{j(\log T)/(KE\log y) \leq i \leq (j+1)(\log T)/(KE\log y) \}$ for some distinct {\em even} $j$, such that
$$ \{z + iE/\log T : i \in \mathcal{I}_{k}\} \subseteq \mathcal{H}^{*}_{\text{even}} \;\;\; \text{and} \;\;\; \#\mathcal{I}_{k} \geq (\log T)/(2KE\log y) \;\;\; \forall 1 \leq k \leq \log y . $$

The point of this manoeuvre is that we have good information about $\p(\max_{i \in \mathcal{I}_{k}} Z(h_{i}) > \sqrt{2(\log\log T - \log\log y)})$ for each set $\mathcal{I}_{k}$, in view of Lemma 3, and moreover the different sets $\mathcal{I}_{k}$ are sufficiently separated that the random variables $Z(h_i)$ corresponding to different sets are ``almost independent\footnote{This is the basic reason for introducing the parameter $y$ in the first place: we will use the blocks $(\mathcal{I}_{k})_{1 \leq k \leq \log y}$ to convert the probability lower bound in Lemma 3, which is fairly large but still $o(1)$, into an overall lower bound $1-o(1)$. There are general concentration inequalities for suprema of Gaussian processes that could also be used for this, but would give a weaker result in Proposition 2.} ''. (Note that $|h_{i}-h_{j}| \geq 1/(K\log y)$ if $h_{i},h_{j}$ correspond to different blocks $\mathcal{I}_{k}$, and recall that the covariances (\ref{correst}) decay rapidly between distant points $h_i , h_j$.)

More precisely, let us write $\max_{k; i \in \mathcal{I}_{k}} Z(h_i)$ to mean $\max_{k=1}^{[\log y]} \max_{i \in \mathcal{I}_{k}} Z(h_{i})$. Then we obviously have that $\p(\max_{h \in \mathcal{H}^{*}} Z(h) > \sqrt{2(\log\log T - \log\log y)})$ is
\begin{eqnarray}
& \geq & \p(\max_{k; i \in \mathcal{I}_{k}} Z(h_{i}) > \sqrt{2(\log\log T - \log\log y)}) \nonumber \\
& = & 1 - \prod_{k=1}^{[\log y]} \p(\max_{i \in \mathcal{I}_{k}} Z(h_{i}) \leq \sqrt{2(\log\log T - \log\log y)}) \nonumber \\
&& + \prod_{k=1}^{[\log y]} \p(\max_{i \in \mathcal{I}_{k}} Z(h_{i}) \leq \sqrt{2(\log\log T - \log\log y)}) - \p(\max_{k; i \in \mathcal{I}_{k}} Z(h_{i}) \leq \sqrt{2(\log\log T - \log\log y)}) , \nonumber
\end{eqnarray}
and in view of Lemma 3 this is
\begin{eqnarray}
& \geq & 1 - \left(1 - \frac{c\sqrt{\log E}}{E (\log\log T)^{3/2}}\right)^{\log y} \nonumber \\
&& + \prod_{k=1}^{[\log y]} \p(\max_{i \in \mathcal{I}_{k}} Z(h_{i}) \leq \sqrt{2(\log\log T - \log\log y)}) - \p(\max_{k; i \in \mathcal{I}_{k}} Z(h_{i}) \leq \sqrt{2(\log\log T - \log\log y)}), \nonumber
\end{eqnarray}
where $c > 0$ is a small absolute constant. In particular, provided that
$$ \frac{\sqrt{\log E} \log y}{E (\log\log T)^{3/2}} \rightarrow \infty \;\;\; \text{as} \; T \rightarrow \infty $$
then the second term here is $o(1)$.

To estimate the difference in the final line we can use the following result, which is one of a family of {\em normal comparison inequalities} that bound a difference between multivariate normal probabilities in terms of differences of the covariance matrices of the relevant random variables. This particular result is due to Li and Shao~\cite{lishao}.
\begin{comparison}[Li and Shao, 2002]
Let $(X_{1},X_{2},...,X_{n})$ and $(W_{1},...,W_{n})$ each be a vector of mean zero, variance one, jointly normal random variables, and write $r_{i,j}^{(1)} = \E X_{i}X_{j}$ and $r_{i,j}^{(0)} = \E W_{i}W_{j}$. Let $u_{1},...,u_{n}$ be any real numbers. Then
\begin{eqnarray}
&& |\p(X_{j} \leq u_{j} \; \forall 1 \leq j \leq n) - \p(W_{j} \leq u_{j} \; \forall 1 \leq j \leq n)| \nonumber \\
& \leq & \frac{1}{2\pi} \sum_{1 \leq i < j \leq n} |\arcsin(r_{i,j}^{(1)}) - \arcsin(r_{i,j}^{(0)})| e^{-(u_{i}^{2} + u_{j}^{2})/(2(1+\max\{|r_{i,j}^{(1)}| , |r_{i,j}^{(0)}|\}))} . \nonumber
\end{eqnarray}
\end{comparison}
The difference that we want to bound is of the kind treated by Comparison Inequality 1, where $u_{j} = \sqrt{2(\log\log T - \log\log y)}$ for all $j$, the $W_{j}$ are simply the random variables $(Z(h_{i}))_{1 \leq k \leq \log y, \; i \in \mathcal{I}_{k}}$, and the $X_{j}$ are the same but with the covariances $\E Z(h_{i_1})Z(h_{i_2})$ replaced by zero when $i_1 , i_2$ do not belong to the same block $\mathcal{I}_{k}$. (Since our random variables are jointly normal, this is equivalent to saying that the $X_{j}$ are the same as the $(Z(h_{i}))$ except that they are independent in different blocks $\mathcal{I}_{k}$, hence the probability factors as a product over $k$.) Thus we find
\begin{eqnarray}
&& | \prod_{k=1}^{[\log y]} \p(\max_{i \in \mathcal{I}_{k}} Z(h_{i}) \leq \sqrt{2(\log\log T - \log\log y)}) - \p(\max_{k; i \in \mathcal{I}_{k}} Z(h_{i}) \leq \sqrt{2(\log\log T - \log\log y)}) | \nonumber \\
& \leq & \frac{1}{2\pi} \sum_{1 \leq k < l \leq \log y} \sum_{i \in \mathcal{I}_{k}} \sum_{j \in \mathcal{I}_{l}} |\arcsin( \E Z(h_{i})Z(h_{j}) )| e^{-2(\log\log T - \log\log y)/(1+ |\E Z(h_{i})Z(h_{j})|) } \nonumber \\
& \ll & \frac{\log^{2}y}{\log^{2}T} \sum_{1 \leq k < l \leq \log y} \sum_{i \in \mathcal{I}_{k}} \sum_{j \in \mathcal{I}_{l}} \frac{1}{|h_i - h_j| \log y \log\log T} , \nonumber
\end{eqnarray}
where the final line uses the estimate $(\ref{correst})$ for the correlations, and we note in particular that $|\E Z(h_{i})Z(h_{j})| \ll 1/\log\log T$ provided $|h_{i}-h_{j}| \gg 1/\log y$, as is the case when $i \in \mathcal{I}_{k}, j \in \mathcal{I}_{l}$ by construction of the blocks $\mathcal{I}_{k}$. Since $\#\mathcal{I}_{k} \leq (\log T)/(E\log y)$ for all $k$, and $|h_{i} - h_{j}| \gg |l-k|/\log y$ for all $i \in \mathcal{I}_{k}, j \in \mathcal{I}_{l}$, the above is
$$ \ll \frac{1}{E^{2}} \sum_{1 \leq k < l \leq \log y} \frac{1}{(l-k) \log\log T} \ll \frac{\log y \log\log y}{E^{2} \log\log T} . $$

Finally, choosing $E=\sqrt{\log\log T} (\log\log\log T)^{2}$ and $\log y = (\log\log T)^{2} (\log\log\log T)^{2}$, say, we have
$$ \frac{\sqrt{\log E} \log y}{E (\log\log T)^{3/2}} \rightarrow \infty \;\;\; \text{and} \;\;\; \frac{\log y \log\log y}{E^{2} \log\log T} \rightarrow 0 \;\;\; \text{as} \; T \rightarrow \infty , $$
and therefore we have
$$ \p(\max_{h \in \mathcal{H}^{*}} Z(h) > \sqrt{2(\log\log T - \log\log y)}) \geq 1 - o(1) . $$
In view of Lemma 2, this suffices to complete the proof of Proposition 2.
\begin{flushright}
Q.E.D.
\end{flushright}

\appendix

\section{Covariance calculations}
In this appendix we perform some variance and covariance calculations, that are necessary for the probabilistic arguments in $\S 3$ but are really just estimates for various sums over primes.

For any fixed $2 \leq P \leq Q \leq T$, let us write
$$  X_{P,Q}(h) := \Re(\sum_{P \leq p \leq Q} \frac{U_{p}}{p^{1/2+ih}} \frac{\log(T/p)}{\log T} ), \;\;\;\;\; 0 \leq h \leq 2\pi . $$
Then we have
\begin{eqnarray}
\E X_{P,Q}(h)^{2} & = & \sum_{P \leq p_{1}, p_{2} \leq Q} \frac{\E \Re(U_{p_{1}} p_{1}^{-ih}) \Re(U_{p_{2}} p_{2}^{-ih}) }{p_{1}^{1/2} p_{2}^{1/2}} \frac{\log(T/p_{1}) \log(T/p_{2})}{\log^{2}T} \nonumber \\
& = & \sum_{P \leq p_{1}, p_{2} \leq Q} \frac{\E (1/2)(U_{p_{1}}p_{1}^{-ih} + \overline{U_{p_{1}}}p_{1}^{ih}) (1/2)(U_{p_{2}}p_{2}^{-ih} + \overline{U_{p_{2}}}p_{2}^{ih}) }{p_{1}^{1/2} p_{2}^{1/2}} \frac{\log(T/p_{1}) \log(T/p_{2})}{\log^{2}T} \nonumber \\
& = & \frac{1}{2} \sum_{P \leq p \leq Q} \frac{1}{p} \frac{\log^{2}(T/p)}{\log^{2}T} , \nonumber
\end{eqnarray}
since $\E U_{p_{1}} U_{p_{2}} = 0$ for all $p_{1},p_{2}$, and $\E U_{p_{1}} \overline{U_{p_{2}}} = 0$ unless $p_{1} = p_{2}$, in which case $\E |U_{p}|^{2} = 1$. It is a standard fact (see e.g. Theorem 2.7 of Montgomery and Vaughan~\cite{mv}) that
$$ \sum_{p \leq x} \frac{1}{p} = \log\log x + b + O\left(\frac{1}{\log x}\right) , \;\;\;\;\; \sum_{p \leq x} \frac{\log p}{p} \ll \log x , \;\;\;\;\; \sum_{p \leq x} \frac{\log^{2}p}{p} \ll \log^{2}x , $$
for a certain constant $b$, and therefore we have
\begin{eqnarray}
\E X_{P,Q}(h)^{2} & = & \frac{1}{2} \left(\sum_{P \leq p \leq Q} \frac{1}{p} - \frac{2}{\log T} \sum_{P \leq p \leq Q} \frac{\log p}{p} + \frac{1}{\log^{2}T} \sum_{P \leq p \leq Q} \frac{\log^{2}p}{p} \right) \nonumber \\
& = & (1/2)(\log\log Q - \log\log P + O(1)) . \nonumber
\end{eqnarray}

Turning to covariances, the same calculations as above show that
\begin{eqnarray}
\E X_{P,Q}(h_{1})X_{P,Q}(h_{2}) & = & \sum_{P \leq p_{1}, p_{2} \leq Q} \frac{\E \Re(U_{p_{1}} p_{1}^{-ih_{1}}) \Re(U_{p_{2}} p_{2}^{-ih_{2}}) }{p_{1}^{1/2} p_{2}^{1/2}} \frac{\log(T/p_{1}) \log(T/p_{2})}{\log^{2}T} \nonumber \\
& = & \frac{1}{2} \sum_{P \leq p \leq Q} \frac{\cos((h_{1}-h_{2})\log p)}{p} \left(1 - \frac{2\log p}{\log T} + \frac{\log^{2}p}{\log^{2}T} \right) . \nonumber
\end{eqnarray}
More explicitly, by a strong form of the prime number theorem (see e.g. Theorem 6.9 of Montgomery and Vaughan~\cite{mv}) we have
$$ \pi(z) := \#\{p \leq z : p \; \text{prime}\} = \int_{2}^{z} \frac{du}{\log u} + O(ze^{-d\sqrt{\log z}}), \;\;\;\;\; z \geq 2, $$
where $d > 0$ is an absolute constant. Therefore, for any $\alpha \neq 0$ we find $\sum_{P \leq p \leq Q} \cos(\alpha\log p)/p$ is
\begin{eqnarray}
\int_{P}^{Q} \frac{\cos(\alpha \log u)}{u} d\pi(u) & = & \int_{P}^{Q} \frac{\cos(\alpha \log u)}{u \log u} du + O((1+|\alpha|)e^{-d\sqrt{\log P}}) \nonumber \\
& = & \int_{\alpha\log P}^{\alpha\log Q} \frac{\cos v}{v}dv + O((1+|\alpha|)e^{-d\sqrt{\log P}}) \nonumber \\
& = & \left\{ \begin{array}{ll}
     \log\log Q - \log\log P + O(1) & \text{if} \; |\alpha\log Q| \leq 1   \\
     \log(1/|\alpha \log P|) + O(1) & \text{if} \; \frac{1}{\log Q} < |\alpha| \leq \frac{1}{\log P}   \\
     O(1/(|\alpha \log P|) + (1+|\alpha|)e^{-d\sqrt{\log P}}) & \text{otherwise,}
\end{array} \right. \nonumber
\end{eqnarray}
the final line using the estimate $\cos v = 1 + O(v^{2})$ when $|v| \leq 1$, and integration by parts on the rest of the range of integration. Similar calculations show that
$$ \sum_{P \leq p \leq Q} \frac{\cos(\alpha \log p) \log p}{p \log T}, \; \sum_{P \leq p \leq Q} \frac{\cos(\alpha \log p) \log^{2}p}{p \log^{2}T}  = O(1/(1+|\alpha \log T|) + (1+|\alpha|)e^{-d\sqrt{\log P}}) , $$
and therefore we have
\begin{equation}
\E X_{P,Q}(h_1)X_{P,Q}(h_2) = \left\{ \begin{array}{ll}
     (1/2)(\log\log Q - \log\log P + O(1)) & \text{if} \; |h_{1}-h_{2}| \leq \frac{1}{\log Q}   \\
     (1/2)(\log(1/|h_{1}-h_{2}|) - \log\log P + O(1)) & \text{if} \; \frac{1}{\log Q} < |h_{1}-h_{2}| \leq \frac{1}{\log P}   \\
     O(1/(|h_{1}-h_{2}| \log P)) & \text{if} \; \frac{1}{\log P} < |h_{1}-h_{2}| \leq 2\pi .
\end{array} \right. \nonumber
\end{equation}

Finally we make a few more qualitative observations about $\E X_{P,Q}(h_1)X_{P,Q}(h_2)$. Firstly, the foregoing calculations show that $\E X_{P,Q}(h)^{2}$ doesn't depend on $h$, and that $\E X_{P,Q}(h_1)X_{P,Q}(h_2)$ is a function of $|h_{1}-h_{2}|$ (in other words the random variables $X_{P,Q}(h)$ are {\em stationary}). Secondly, there exists an absolute constant $K > 0$ such that
$$ \E X_{P,Q}(h_1)X_{P,Q}(h_2) \geq 0 \;\;\;\;\; \text{if} \; |h_{1}-h_{2}| \leq \frac{1}{K \log P} . $$
And thirdly, if we write $r_{P,Q}(h) := \E X_{P,Q}(h_1)X_{P,Q}(h_1 + h)$, (which doesn't depend on $h_{1}$, because of stationarity), and if $0 < \delta \leq h \leq 2\pi$ and $\sqrt{\log P} \gg \log\log Q$, then we have
\begin{eqnarray}
r_{P,Q}(h) - r_{P,Q}(h+\delta) & = & \frac{1}{2} \int_{h\log P}^{(h+\delta)\log P} \frac{\cos v}{v} dv - \frac{1}{2} \int_{h\log Q}^{(h+\delta)\log Q} \frac{\cos v}{v} dv + O(\frac{1}{h\log T} + e^{-d\sqrt{\log P}}) \nonumber \\
& = & \frac{1}{2} \int_{h\log P}^{(h+\delta)\log P} \frac{\cos v}{v} dv + O(\frac{1}{h\log Q}) , \nonumber
\end{eqnarray}
using integration by parts. In particular, if $h+\delta \leq 1/\log P$ (say) then the integral here is $\gg \int_{h\log P}^{(h+\delta)\log P} \frac{1}{v} dv \gg \delta/h $, and so if $\delta \geq K/\log Q$ then
$$ r_{P,Q}(h) - r_{P,Q}(h+\delta) > 0. $$

\section{A multivariate central limit theorem}
In this appendix we formulate a version of the central limit theorem that justifies Lemma 2, above, in which the random variables $Y(h), h \in \mathcal{H}^{*}$ were replaced by normal random variables $Z(h)$ with the same means and covariances.

In fact, we will sketch a proof of the following theorem.
\begin{clt}[Specialised from Theorem 2.1 of Reinert and R\"{o}llin~\cite{reinertrollin}]
Suppose that $n \geq 1$, and that $\mathcal{H}$ is a finite non-empty set. Suppose that for each $1 \leq i \leq n$ and $h \in \mathcal{H}$ we are given a deterministic coefficient $c(i,h) \in \C$. Finally, suppose that $(V_{i})_{1 \leq i \leq n}$ is a sequence of independent, mean zero, complex valued random variables, and let $Y = (Y_{h})_{h \in \mathcal{H}}$ be the $\#\mathcal{H}$-dimensional random vector with components
$$ Y_{h} := \Re\left(\sum_{i=1}^{n} c(i,h) V_{i} \right). $$

If $Z = (Z_{h})_{h \in \mathcal{H}}$ is a multivariate normal random vector with the same mean vector and covariance matrix as $Y$, then for any $u \in \R$ and any small $\delta > 0$ we have
\begin{eqnarray}
\p(\max_{h \in \mathcal{H}} Y_{h} \leq u) & \leq & \p(\max_{h \in \mathcal{H}} Z_{h} \leq u + \delta) + \nonumber \\
&& + O\left(\frac{1}{\delta^{2}} \sum_{g,h \in \mathcal{H}} \sqrt{\sum_{i=1}^{n} |c(i,g)|^{2}|c(i,h)|^{2} \E|V_{i}|^{4}} + \frac{1}{\delta^{3}} \sum_{i=1}^{n} \E|V_{i}|^{3} \left(\sum_{h \in \mathcal{H}} |c(i,h)| \right)^{3} \right). \nonumber
\end{eqnarray}
There is an exactly similar lower bound for $\p(\max_{h \in \mathcal{H}} Y_{h} \leq u)$, in which $\p(\max_{h \in \mathcal{H}} Z_{h} \leq u + \delta)$ is replaced by $\p(\max_{h \in \mathcal{H}} Z_{h} \leq u - \delta)$.
\end{clt}

Before discussing its proof, let us see how Central Limit Theorem 1 is applicable to the random variables $(Y(h))_{h \in \mathcal{H}^{*}}$ considered in $\S 3$. Recall that, by definition,
$$ Y(h) := \frac{\Re \sum_{y \leq p \leq T} \frac{U_{p}}{p^{1/2+ih}} \frac{\log(T/p)}{\log T} }{\sqrt{\frac{1}{2} \sum_{y \leq p \leq T} \frac{1}{p} \frac{\log^{2}(T/p)}{\log^{2}T}}} , $$
so we are in the setting of Central Limit Theorem 1 with the indices $i$ replaced by primes $y \leq p \leq T$, and $V_{p} = U_{p}$, and
$$ c(p,h) = \frac{\frac{1}{p^{1/2+ih}} \frac{\log(T/p)}{\log T}}{\sqrt{\frac{1}{2} \sum_{y \leq p \leq T} \frac{1}{p} \frac{\log^{2}(T/p)}{\log^{2}T}}} = \frac{\frac{1}{p^{1/2+ih}} \frac{\log(T/p)}{\log T}}{\sqrt{(1/2)(\log\log T - \log\log y + O(1))}} . $$
Now Central Limit Theorem 1 implies that $\p(\max_{h \in \mathcal{H}^{*}} Y(h) \leq u)$ is at most
$$ \p(\max_{h \in \mathcal{H}^{*}} Z(h) \leq u + \delta) + O\left( \frac{(\#\mathcal{H}^{*})^{2}}{\delta^{2}} \sqrt{\sum_{y \leq p \leq T} \frac{1}{p^{2} (\log\log T)^{2}}} + \frac{(\#\mathcal{H}^{*})^{3}}{\delta^{3}} \sum_{y \leq p \leq T} \frac{1}{p^{3/2} (\log\log T)^{3/2}} \right) . $$
Provided that $y \geq ((\#\mathcal{H}^{*})/\delta)^{6}$, say, the ``big Oh'' term here is $o(1)$ as $T \rightarrow \infty$. Choosing $u = (\log\log T - 2\log\log\log T - (\log\log\log T)^{3/4})/(\sqrt{(1/2)(\log\log T - \log\log y + O(1))})$, and $\delta = 1/\sqrt{\log\log T}$, we conclude that
\begin{eqnarray}
&& \p(\max_{h \in \mathcal{H}^{*}} Y(h) > \frac{\log\log T - 2\log\log\log T - (\log\log\log T)^{3/4}}{\sqrt{(1/2)(\log\log T - \log\log y + O(1))}}) \nonumber \\
& \geq & \p(\max_{h \in \mathcal{H}^{*}} Z(h) > \frac{\log\log T - 2\log\log\log T - (\log\log\log T)^{3/4}}{\sqrt{(1/2)(\log\log T - \log\log y + O(1))}} + \frac{1}{\sqrt{\log\log T}}) + o(1), \nonumber
\end{eqnarray}
provided that $y \geq (\log\log T)^{3}(\#\mathcal{H}^{*})^{6}$. Since $\#\mathcal{H}^{*} \ll (\log T)/E$, this condition is certainly satisfied if $y \geq \log^{7}T$, say, and so Lemma 2 follows.

\vspace{12pt}
We will deduce Central Limit Theorem 1 from a much more general normal approximation theorem of Reinert and R\"{o}llin~\cite{reinertrollin}, which they prove using Stein's method of exchangeable pairs. If one only wants a result like Central Limit Theorem 1 this could probably be deduced from many other existing results as well, but since Reinert and R\"{o}llin's result is neat and powerful, and automatically supplies explicit ``big Oh'' terms, we will work from there. See also the second appendix in the author's paper~\cite{harpergp}, where Reinert and R\"{o}llin's result was applied to a very similar problem.

To apply Theorem 2.1 of Reinert and R\"{o}llin~\cite{reinertrollin}, we first need to construct a random vector $Y' = (Y'_{h})_{h \in \mathcal{H}}$ such that the pair $(Y,Y')$ is {\em exchangeable}, i.e. such that $(Y,Y')$ has the same law as $(Y',Y)$. In fact there is a standard way to do this: let $I$ be a random variable independent of everything else, having the discrete uniform distribution on the set $\{1,2,...,n\}$, and let $(V_{i}')_{1 \leq i \leq n}$ be independent random variables having the same distribution as the $V_{i}$; then conditional on the event $I=i$, define
$$ Y'_{h} := Y_{h} - \Re\left(c(i,h) V_{i} \right) + \Re\left(c(i,h) V'_{i} \right) , \;\;\;\;\; h \in \mathcal{H}. $$
Since there is perfect symmetry between the roles of $V_{i}$ and $V'_{i}$, the reader may readily convince themselves that $(Y,Y')$ form an exchangeable pair.

Next we shall perform a few conditional expectation calculations we will need, beginning with the calculation of $\E(Y'-Y |Y) = \E((Y'_{h} - Y_{h})_{h \in \mathcal{H}} | (Y_{h})_{h \in \mathcal{H}})$. Since $I$ is independent of $Y$, and is distributed uniformly on $\{1,2,...,n\}$, we have
\begin{eqnarray}
\E(Y'-Y |Y) & = & \sum_{i=1}^{n} \frac{1}{n} \E(Y'-Y |Y, I=i) \nonumber \\
& = & \sum_{i=1}^{n} \frac{1}{n} \E((- \Re\left(c(i,h) V_{i} \right) + \Re\left(c(i,h) V'_{i} \right))_{h \in \mathcal{H}} | (Y_{h})_{h \in \mathcal{H}}, I=i) \nonumber \\
& = & \sum_{i=1}^{n} \frac{1}{n} \E((- \Re\left(c(i,h) V_{i} \right))_{h \in \mathcal{H}} | (Y_{h})_{h \in \mathcal{H}}, I=i) , \nonumber
\end{eqnarray}
the final line using the fact that $V_{i}'$ is independent of $Y$ and $I$, and has mean zero. Now we observe that $(- \Re\left(c(i,h) V_{i} \right))_{h \in \mathcal{H}}$ is independent of $I$, and so we have
\begin{eqnarray}
\E(Y'-Y |Y) = \sum_{i=1}^{n} \frac{1}{n} \E((- \Re\left(c(i,h) V_{i} \right))_{h \in \mathcal{H}} | (Y_{h})_{h \in \mathcal{H}}) & = & \frac{1}{n} \E(\sum_{i=1}^{n} (- \Re\left(c(i,h) V_{i} \right))_{h \in \mathcal{H}} | (Y_{h})_{h \in \mathcal{H}}) \nonumber \\
& = & \frac{-1}{n} \E((Y_{h})_{h \in \mathcal{H}} |(Y_{h})_{h \in \mathcal{H}}) , \nonumber
\end{eqnarray}
where the second equality uses linearity of expectation, and the third equality uses the definition of $Y_{h}$. We conclude that
$$ \E(Y'-Y |Y) = \frac{-1}{n} Y . $$
In a similar way, for any $g,h \in \mathcal{H}$ we have
\begin{eqnarray}
&& \E((Y'_{g}-Y_{g})(Y'_{h}-Y_{h}) | Y) \nonumber \\
& = & \frac{1}{n} \E(\sum_{i=1}^{n} (- \Re\left(c(i,g) V_{i} \right) + \Re\left(c(i,g) V_{i}' \right)) (- \Re\left(c(i,h) V_{i} \right) + \Re\left(c(i,h) V_{i}' \right)) | Y) \nonumber \\
& = & \frac{1}{n} \E(\sum_{i=1}^{n} \Re\left(c(i,g) V_{i} \right) \Re\left(c(i,h) V_{i} \right) | Y) + \frac{1}{n} \E(\sum_{i=1}^{n} \Re\left(c(i,g) V_{i}' \right) \Re\left(c(i,h) V_{i}' \right)) , \nonumber
\end{eqnarray}
the second equality using the fact that each $V_{i}'$ is independent of $V_{i}$ and $Y$ (and has mean zero). Also, for any $f,g,h \in \mathcal{H}$ we have
$$ \E|(Y'_{f}-Y_{f})(Y'_{g}-Y_{g})(Y'_{h}-Y_{h})| = \sum_{i=1}^{n} \frac{1}{n} \E|\Re c(i,f)(V_{i}'-V_{i})| |\Re c(i,g)(V_{i}'-V_{i})| |\Re c(i,h)(V_{i}'-V_{i})| . $$

Now Theorem 2.1 of Reinert and R\"{o}llin~\cite{reinertrollin} asserts that, if $t : \R^{\#\mathcal{H}} \rightarrow \R$ is any three times differentiable function,
\begin{eqnarray}
|\E t(Y) - \E t(Z)| & \leq & \frac{1}{4} \left(\sup_{g,h \in \mathcal{H}} || \frac{\partial^{2}}{\partial x_{g} \partial x_{h}} t ||_{\infty} \right) \sum_{g,h \in \mathcal{H}} n \sqrt{\text{Var}\left( \E((Y'_{g}-Y_{g})(Y'_{h}-Y_{h}) | Y) \right)} + \nonumber \\
&& + \frac{1}{12}  \left(\sup_{f,g,h \in \mathcal{H}} || \frac{\partial^{3}}{\partial x_{f} \partial x_{g} \partial x_{h}} t ||_{\infty} \right) \sum_{f,g,h \in \mathcal{H}} n \E|(Y'_{f}-Y_{f})(Y'_{g}-Y_{g})(Y'_{h}-Y_{h})|  . \nonumber
\end{eqnarray}
Here the factor $n$ is the reciprocal of the $1/n$ arising in the condition $ \E(Y'-Y |Y) = \frac{-1}{n} Y $.

To understand the terms in this bound, we note that the second sum in our expression for $\E((Y'_{g}-Y_{g})(Y'_{h}-Y_{h}) | Y)$, above, is deterministic (it is just an expectation, rather than a conditional expectation), so can be ignored when computing the variance of $\E((Y'_{g}-Y_{g})(Y'_{h}-Y_{h}) | Y)$. Thus we have
\begin{eqnarray}
\text{Var}\left( \E((Y'_{g}-Y_{g})(Y'_{h}-Y_{h}) | Y) \right) & = & \text{Var}\left( \frac{1}{n} \E(\sum_{i=1}^{n} \Re\left(c(i,g) V_{i} \right) \Re\left(c(i,h) V_{i} \right) | Y) \right) \nonumber \\
& = & \frac{1}{n^{2}} \text{Var}\left( \E(\sum_{i=1}^{n} \Re\left(c(i,g) V_{i} \right) \Re\left(c(i,h) V_{i} \right) | Y) \right) \nonumber \\
& \leq & \frac{1}{n^{2}} \text{Var}\left( \sum_{i=1}^{n} \Re\left(c(i,g) V_{i} \right) \Re\left(c(i,h) V_{i} \right) \right) , \nonumber
\end{eqnarray}
where the final line uses the fact that conditioning reduces variance. At this point, since the $V_{i}$ are independent we conclude that $\text{Var}\left( \E((Y'_{g}-Y_{g})(Y'_{h}-Y_{h}) | Y) \right)$ is
$$ \leq \frac{1}{n^{2}} \sum_{i=1}^{n} \text{Var} \Re\left(c(i,g) V_{i} \right) \Re\left(c(i,h) V_{i} \right) \leq \frac{1}{n^{2}} \sum_{i=1}^{n} \E|c(i,g) V_{i} c(i,h) V_{i}|^{2} = \frac{1}{n^{2}} \sum_{i=1}^{n} |c(i,g)|^{2} |c(i,h)|^{2} \E |V_{i}|^{4} . $$
Rather more straightforwardly, our expression for $\E|(Y'_{f}-Y_{f})(Y'_{g}-Y_{g})(Y'_{h}-Y_{h})|$, above, together with the fact that $\E |V_{i}'|^{3} = \E |V_{i}|^{3}$, imply that
$$ \E|(Y'_{f}-Y_{f})(Y'_{g}-Y_{g})(Y'_{h}-Y_{h})| \ll \sum_{i=1}^{n} \frac{1}{n} |c(i,f)||c(i,g)||c(i,h)| \E|V_{i}|^{3} . $$

Finally we can set $t((x_{h})_{h \in \mathcal{H}}) := \prod_{h \in \mathcal{H}} s(x_{h})$, where $s : \R \rightarrow [0,1]$ is any three times differentiable function such that
$$ s(x) = \left\{ \begin{array}{ll}
     1 & \text{if} \; x \leq u   \\
     0 & \text{if} \; x \geq u+\delta .
\end{array} \right. $$
We can find such $s$ with derivatives satisfying $|s^{(r)}(x)| = O(\delta^{-r})$, $0 \leq r \leq 3$, in which case we will have
$$ \sup_{g,h \in \mathcal{H}} || \frac{\partial^{2}}{\partial x_{g} \partial x_{h}} t ||_{\infty} \ll \delta^{-2}, \;\;\;\;\; \sup_{f,g,h \in \mathcal{H}} || \frac{\partial^{3}}{\partial x_{f} \partial x_{g} \partial x_{h}} t ||_{\infty} \ll \delta^{-3} . $$
Then we see
\begin{eqnarray}
&& \p(\max_{h \in \mathcal{H}} Y_{h} \leq u) - \p(\max_{h \in \mathcal{H}} Z_{h} \leq u + \delta) \nonumber \\
& \leq & \E t(Y) - \E t(Z) \nonumber \\
& \ll & \frac{1}{\delta^{2}} \sum_{g,h \in \mathcal{H}} \sqrt{\sum_{i=1}^{n} |c(i,g)|^{2}|c(i,h)|^{2} \E|V_{i}|^{4}} + \frac{1}{\delta^{3}} \sum_{i=1}^{n} \E|V_{i}|^{3} \left(\sum_{h \in \mathcal{H}} |c(i,h)| \right)^{3} , \nonumber
\end{eqnarray}
which is the upper bound claimed in Central Limit Theorem 1. The lower bound follows by instead choosing $s(x)$ to be 1 if $x \leq u-\delta$, and 0 if $x \geq u$.
\begin{flushright}
Q.E.D.
\end{flushright}

\vspace{12pt}
\noindent {\em Acknowledgements.}  The author would like to thank Louis-Pierre Arguin for drawing his attention to Fyodorov and Keating's work, and for discussions about logarithmically correlated random variables. He would also like to thank Maksym Radziwi\l\l , for suggestions about rigorous results on the maximum of the zeta function, and Yan Fyodorov and Jon Keating for their encouragement and help with the references.


\begin{thebibliography}{99}

\bibitem{diaconisshah} P. Diaconis, M. Shahshahani. On the eigenvalues of random matrices. {\em J. Appl. Probab.}, \textbf{31A} (Studies in applied probability), pp 49-62. 1994

\bibitem{farmergonekhughes} D. W. Farmer, S. M. Gonek, C. P. Hughes. The maximum size of $L$-functions. {\em J. Reine Angew. Math.,} \textbf{609}, pp 215-236. 2007

\bibitem{fyodhiarykeat} Y. V. Fyodorov, G. A. Hiary, J. P. Keating. Freezing Transition, Characteristic Polynomials of Random Matrices, and the Riemann Zeta-Function. {\em Phys. Rev. Lett.}, \textbf{108}, 170601 (5pp). 2012

\bibitem{fyodkeat} Y. V. Fyodorov, J. P. Keating. Freezing Transitions and Extreme Values: Random Matrix Theory, $\zeta(1/2+it)$, and Disordered Landscapes. To appear in {\em Phil. Trans. R. Soc. A}. Preprint available online at \url{http://arxiv.org/abs/1211.6063}

\bibitem{gonekhugheskeating} S. M. Gonek, C. P. Hughes, J. P. Keating. A hybrid Euler-Hadamard product for the Riemann zeta function. {\em Duke Math. J.}, \textbf{136}, no. 3, pp 507-549. 2007

\bibitem{harpergp} A. J. Harper. Bounds on the suprema of Gaussian processes, and omega results for the sum of a random multiplicative function. {\em Ann. Appl. Probab.}, \textbf{23}, no. 2, pp 584-616. 2013

\bibitem{ivic} A. Ivi\'{c}. {\em The Riemann Zeta-Function: Theory and Applications.} Dover republished edition, published by Dover Publications, Inc.. 2003

\bibitem{lishao} W. Li, Q.-M. Shao. A normal comparison inequality and its applications. {\em Probab. Theory Relat. Fields}, \textbf{122}, pp 494-508. 2002

\bibitem{mv} H. L. Montgomery, R. C. Vaughan. {\em Multiplicative Number Theory I: Classical Theory.} First edition, published by Cambridge University Press. 2007

\bibitem{reinertrollin} G. Reinert, A. R\"{o}llin. Multivariate Normal Approximation with Stein's Method of Exchangeable Pairs under a General Linearity Condition. {\em Ann. Probab.}, \textbf{37}, no. 6, pp 2150-2173. 2009

\bibitem{soundmoments} K. Soundararajan. Moments of the Riemann zeta function. {\em Ann. Math.}, \textbf{170}, pp 981-993. 2009

\bibitem{talagrand} M. Talagrand. The Missing factor in Hoeffding's inequalities. {\em Ann. Inst. Henri Poincar\'{e} - Probabilit\'{e}s et Statistiques}, \textbf{31}, no. 4, pp 689-702. 1995

\bibitem{zeitouni} O. Zeitouni. Branching random walks and Gaussian fields. {\em Notes for Lectures}, available online at \url{http://www.wisdom.weizmann.ac.il/~zeitouni}



\end{thebibliography}
\end{document}